\documentclass[12pt]{amsart}

\usepackage{graphicx}
\usepackage{stmaryrd}
\usepackage{amssymb}

\usepackage{amsmath,amscd}

\usepackage[T1]{fontenc}
\usepackage{yfonts}

\renewcommand{\iff}{\mbox{if{}f}}

\newcommand{\spc}{(\omega+1)^\omega}
\newcommand{\tr}{(\omega+1)^{<\omega}}

\newcommand{\har}{\!\!\upharpoonright\!\!}
\renewcommand{\int}{\mbox{int}}
\newcommand{\diam}{\mbox{diam}}
\newcommand{\dom}{\mbox{dom}}
\newcommand{\cl}{\mbox{cl}}
\newcommand{\otr}{\omega^{<\omega}}
\renewcommand{\oplus}{(\omega+1)}
\newcommand{\Bor}{\mbox{Bor}}

\newtheorem{theorem}{Theorem}
\newtheorem{lemma}{Lemma}
\newtheorem{corollary}{Corollary}

\theoremstyle{definition}

\author{Marcin Sabok} 

\address{Mathematical Institute,
  Wroc\l aw University, pl. Grunwaldzki $2\slash 4$,
  $50$-$384$ Wroc\l aw, Poland }

\email{sabok@math.uni.wroc.pl}

\title{A dichotomy for Borel functions}

\begin{document}

\begin{abstract}
  The dichotomy discovered by Solecki in \cite{Sol} states
  that any Baire class 1 function is either
  $\sigma$-continuous or ``includes'' the Pawlikowski
  function $P$. The aim of this paper is to give an argument
  which is simpler than the original proof of Solecki and
  gives a stronger statement: a dichotomy for all Borel
  functions.
\end{abstract}

\maketitle

\section{Introduction}

An old question of Lusin asked whether there exists a Borel
function which cannot be decomposed into countably many
continuous functions. By now several examples have been
given, by Keldi\v{s}, Adyan and Novikov among others. A
particularly simple example, the function
$P:\spc\rightarrow\omega^\omega$, has been found by
Pawlikowski (cf. \cite{CMPS}). By definition,
\begin{displaymath}
  P(x)(n) = \left\{ 
    \begin{array}{ll}
      x(n)+1 & \mbox{if}\quad x(n)<\omega,\\
      0 & \mbox{if}\quad x(n)=\omega.
    \end{array} \right.
\end{displaymath}
It is proved in \cite{CMPS} that if $A\subseteq\spc$ is such
that $P\har A$ is continuous then
$P[A]\subseteq\omega^\omega$ is nowhere dense. Since $P$ is
a surjection, it is not $\sigma$-continuous.

In \cite{Sol} Solecki showed that the above function is, in
a sense, the only such example, at least among Baire class 1
functions (in other words, it is the initial object in a
certain category).

\begin{theorem}[Solecki, \cite{Sol}]\label{pmin}
  For any Baire class 1 function $f:X\rightarrow Y$, where
  $X,Y$ are Polish spaces, either $f$ is $\sigma$-continuous
  or there exist topological embeddings $\varphi$ and $\psi$
  such that the following diagram commutes:
  $$
  \begin{CD}
    \omega^\omega @ >\psi >> Y\\
    @AAP A @AAfA\\
    (\omega+1)^\omega @> \varphi >> X
  \end{CD}
  $$
\end{theorem}

In \cite{Zpl:DSTDF} Zapletal generalized Solecki's dichotomy
to all Borel functions by proving the following theorem.

\begin{theorem}[Zapletal, \cite{Zpl:DSTDF}]
  If $f:X\rightarrow Y$ is a Borel function which is not
  $\sigma$-continuous then there is a compact set
  $C\subseteq X$ such that $f\har C$ is not
  $\sigma$-continuous and of Baire class 1.
\end{theorem}

In this paper we give a new proof of the above dichotomy for
all Borel functions, which is direct, shorter and more
general than the original proof from \cite{Sol}.

\section{Notation}

We say that a Borel function $f:X\rightarrow Y$, where $X,Y$
are Polish spaces, is $\sigma$-continuous if there exist a
countable cover of the space $X=\bigcup_n X_n$ (with
arbitrary sets $X_n$) such that $f\har X_n$ is continuous
for each $n$. It follows from the Kuratowski extension
theorem that we may require that the sets $X_n$ be Borel. If
$f$ is a Borel function which is not $\sigma$-continuous
then the family of sets on which it is $\sigma$-continuous
is a proper $\sigma$-ideal in $X$. We denote this
$\sigma$-ideal by $I_f$.

In a metric space $(X,d)$ for $A,B\subseteq X$ let us denote
by $h(A,B)$ the Hausdorff distance between $A$ and $B$.

The spaces $\spc$ and $\oplus^n$ for $n<\omega$ are endowed
with the product topology of order topologies on $\omega+1$.

\section{The Zapletal's game}

In \cite{Zpl:DSTDF} Zapletal introduced a two-player game,
which turnes out to be very useful in examining
$\sigma$-continuity of Borel functions. Let
$B\subseteq\omega^\omega$ be a Borel set and $f:B\rightarrow
2^\omega$ be a Borel function. Let
$\rho:\omega\rightarrow\omega\times 2^{<\omega}\times\omega$
be a bijection. The game $G_f(B)$ is played by Adam and Eve.
They take turns playing natural numbers. In his $n$-th move,
Adam picks $x_n\in\omega$. In her $n$-th move, Eve chooses
$y_n\in 2$. At the end of the game we have
$x\in\omega^\omega$ and $y\in 2^\omega$ formed by the
numbers picked by Adam and Eve, respectively. Next, $y\in
2^\omega$ is used to define a sequence of partial continuous
functions (with domains of type $G_\delta$ in
$\omega^\omega$) in the following way. For $n<\omega$ let
$f_n$ be a partial function from $\omega^\omega$ to
$2^\omega$ such that for $t\in\omega^\omega$ and $\sigma\in
2^{<\omega}$
$$ f_n(t)\supseteq\sigma\quad\iff\quad\exists
k\in\omega\,\,\, y(\rho(n,\sigma,k))=1$$ and
$\dom(f_n)=\{t\in\omega^\omega: \forall
n<\omega\,\exists!\sigma\in 2^n\,\,
f_n(t)\supseteq\sigma\}$.  Eve wins the game $G_f(B)$ if
$x\not\in B$ or $\exists n\, f(x)=f_n(x).$ Otherwise Adam
wins the game.

It is easy to see that if $f$ is a Borel function then $G_f$
is a Borel game. The key feature of the game $G_f$ is that
it detects $\sigma$-continuity of the function $f$.

\begin{theorem}[Zapletal,\cite{Zpl:DSTDF}]
  For $B\subseteq\omega^\omega$ and $f:B\rightarrow
  2^\omega$ Eve has a winning strategy in the game $G_f(B)$
  if and only if $f$ is $\sigma$-continuous on $B$.
\end{theorem}

Note that if Adam has a winning strategy then the image of
his strategy (treated as a continuous function from
$2^\omega$ to $B$) is a compact set on which $f$ is also not
$\sigma$-continuous. This observation and the Borel
determinacy gives the following corollary.

\begin{corollary}[Zapletal,\cite{Zpl:DSTDF}]\label{cldense}
  If $B$ is a Borel set and $f:B\rightarrow 2^\omega$ is a
  Borel function which is not $\sigma$-continuous then there
  is a compact set $C\subseteq B$ such that $f\har C$ is
  also not $\sigma$-continuous.
\end{corollary}

\section{Proof of the dichotomy}

In the statement of Theorem \ref{pmin} both functions
$\varphi$ and $\psi$ are to be topological embeddings.
However, as we will see below, for the dichotomy it is
enough that they both are injective, $\varphi$ continuous
and $\psi$ open. We are going to prove first this version of
the dichotomy.

\begin{theorem}\label{dichotomy}
  Let $X$ be a Polish space and $f:X\rightarrow 2^\omega$ be
  a Borel function. Then precisely one of the following
  conditions holds:
 \begin{enumerate}
 \item either $f$ is $\sigma$-continuous
 \item\label{fact} or there are an open injection $\psi$ and
   a continuous injection $\varphi$ such that the following
   diagram commutes:
   $$
     \begin{CD}
       \omega^\omega @ >\psi >> 2^\omega\\
       @AAP A @AAf A\\
       (\omega+1)^\omega @> \varphi >> X
     \end{CD}
   $$   
 \end{enumerate}
\end{theorem}
\noindent Notice that compactness of $\spc$ implies that the
$\psi$ above must be a topological embedding.
\begin{proof}
  
  It is straightforward that (\ref{fact}) implies that $f$
  is not $\sigma$-continuous. Let us assume that $f$ is not
  $\sigma$-continuous and prove that (\ref{fact}) holds.  By
  Corollary \ref{cldense} we may assume that $X$ is compact.

  \subsection*{Notation.}

  First we introduce some notation. For a fixed $n$ and
  $0\leq k\leq n$ let $S^n_k$ be the set of points in
  $\oplus^n$ of Cantor-Bendixson rank $\geq n-k$. For each
  $n<\omega$ and $1\leq k\leq n$ let us pick a function
  $\pi^n_k:S^n_k\rightarrow S^n_{k-1}$ such that
  \begin{itemize}
  \item on $S^n_{k-1}$ $\pi^n_k$ is the identity,
  \item if $\tau\in S^n_k\setminus S^n_{k-1}$ then we pick
    one $i\in n$ such that $\tau(i)<\omega$ and $\tau(i)$ is
    maximal such and define
    $$\pi^n_k(\tau)(i)=\omega, \quad
    \pi^n_k(\tau)(j)=\tau(j) \ \ \ \mbox{for}\ j\not= i.$$
  \end{itemize} 
  This definition clearly depends on the choice of the index
  $i$ above. Note, however, that we may pick the functions
  $\pi^n_k$ so that they are coherent, in the sense that for
  $\tau\in\oplus^{n+1}$, unless $\tau(n)$ is the biggest
  finite value of $\tau$, we have
  $\pi^{n+1}_{k+1}(\tau)=\pi^n_k(\tau\har
  n)^\smallfrown\tau(n)$. In particular
  $\pi^{n+1}_{k+1}(\sigma^\smallfrown\omega)=\pi^n_k(\sigma)^\smallfrown\omega$
  for any $\sigma\in\oplus^n$.  The functions $\pi^n_k$ will
  be called projections.

  \begin{lemma}
    For each $n$ and $1\leq k\leq n$ the projection
    $\pi^n_k: S^n_k\rightarrow S^n_{k-1}$ is continuous.
  \end{lemma}
  \begin{proof}
    Note that any point in $S^n_k$ except
    $(\omega,\ldots,\omega)$ ($k$ times $\omega$) has a
    neighborhood in which projection is unambigous and hence
    continuous. But it is easy to see that at the point
    $(\omega,\ldots,\omega)$ any projection is continuous.
  \end{proof}

  For each $n<\omega$ let us also introduce the function
  $r_n:\oplus^n\rightarrow\oplus^{n}$ defined as
  $r_n(\tau^\smallfrown a)=\tau^\smallfrown\omega$.

  To make the above notation more readable we will usually
  drop subscripts and superscripts in $\pi^n_k$ and $r_n$.

  We pick a well-ordering $\leq$ of $\tr$ into type $\omega$
  such that for each point $\tau\in\tr$ all elements of the
  transitive closure of $\tau$ with respect to $\pi$, $r$
  and restrictions (i.e. functions of the form
  $\oplus^n\ni\tau\mapsto\tau\har m\in\oplus^m$ for $m<n$)
  are $\leq\tau$.

  For a set $B\in\Bor(X)\setminus I_f$ let $B^*$ denote the
  set $B$ shrunk by all basic clopens $C$ which have
  $I$-small intersection with $B$.

  \subsection*{Strategy of the construction.}
  
  In order to define functions $\varphi$ and $\psi$, we will
  construct for each $\tau\in\tr$ a clopen set
  $C_\tau\subseteq 2^\omega$ and a compact set
  $X_\tau\subseteq X$ such that if $\sigma\subseteq\tau$
  then $C_\tau\subseteq C_\sigma$ and $X_\tau\subseteq
  X_\sigma$.

  The sets $C_\tau$ will be disjoint, which means that for
  $\tau\not=\tau'$, $|\tau|=|\tau'|$ $C_\tau\cap
  C_{\tau'}=\emptyset$. We will also need $X_\tau\subseteq
  f^{-1}[C_\tau]$ and $\diam(X_\tau)<1\slash |\tau|$.
 
  The construction of the sets $X_\tau, C_\tau$ will be done
  by induction along the ordering $\leq$ on $\tr$.  In fact,
  we will do something more: at each step $n$ if $\tau$ is
  the $n$-th element of $\tr$ we will construct not only a
  compact set $X_\tau$ but also $I_f$-positive Borel sets
  $X^n_\sigma$ for $\sigma\leq\tau$ such that:
  \begin{itemize}
  \item $X_\tau\subseteq
    X^{n-1}_{\tau\upharpoonright(|\tau|-1)}$,
  \item $X^n_\sigma\subseteq X^{n-1}_\sigma$ if
    $\sigma<\tau$,
  \item $X^n_{\sigma^\smallfrown a}\subseteq X^n_\sigma$ if
    $\sigma,\sigma^\smallfrown a<\tau$,
  \item $X^n_\sigma\cap f^{-1}[C_{\sigma^\smallfrown
      a}]=\emptyset$ if $\sigma,\sigma^\smallfrown
    a\leq\tau$,
  \item
    $X^n_{\sigma^\smallfrown\omega}\subseteq\cl(X^n_\sigma)$
    if $\sigma,\sigma^\smallfrown\omega\leq\tau$.
  \end{itemize}
  The set $X^n_\sigma$ is to be understood as the space for
  further construction of sets $X_\rho$ for
  $\rho\supseteq\sigma$ and $\rho>\tau$, as can be seen in
  the first condition above. The last condition, as we will
  see later, will be used to guarantee ``continuity'' of the
  family of sets $X_\tau$. For technical reasons we will
  also make sure that $X^n_\sigma=(X^n_\sigma)^*$.

  We are going to ensure disjointness of $C_\tau$'s by
  satisfying the following conditions:
  \begin{itemize}
  \item $C_{\tau^\smallfrown a}\subseteq C_\tau$,
  \item $C_{\tau^\smallfrown a}\cap C_{\tau^\smallfrown
       b}=\emptyset$ for $a\not=b$.
  \end{itemize}

  The fact that $\diam(X_\tau)<1\slash|\tau|$ will follow
  from the following inductive conditions (recall that
  $\pi(\tau)\leq\tau$ for any $\tau$):
  \begin{itemize}
  \item $\diam(X_\tau)< 3\,\diam(X_{\pi(\tau)})$,
  \item $\diam(X_{\tau^\smallfrown\omega}) <
    1\slash(3^{|\tau|+1}(|\tau|+1))$,
  \end{itemize}
  because iterating projections in $\oplus^n$ stabilizes
  before $n+1$ steps.

  The crucial feature of the sets $X_\tau$ is that this
  family should be ``continuous''. Namely, we will require
  that if $\tau$ and $\pi(\tau)$ occur by the $n$-th step
  then
  \begin{equation}\label{ineq}
    h(X^n_\tau,X^n_{\pi(\tau)})<3^{|\tau|}\,d(\tau,\pi(\tau))
  \end{equation}

  This condition is the most diffucult. To fulfill it we
  will construct yet another kind of objects. Notice first
  that if $h(A,B)<\varepsilon$ for two nonempty sets in $X$
  then there are two finite families (we will refer to them
  as to ``anchors'') $A_i$ and $B_i$ ($i\in I_0$) of subsets
  of $A$ and $B$ respectively such that for any
  $A_i'\subseteq A_i$, $B_i'\subseteq B_i$ still
  $h(\bigcup_i A_i',\bigcup_i B_i')<\varepsilon$. Similarly,
  if $h(A,B)<\varepsilon$ and $C\subseteq A$ then there
  exist a finite family $D_i$ ($i\in I_0$) of subsets of $B$
  such that for any $D_i'\subseteq D_i$ $h(\bigcup_i D_i',
  C)<\varepsilon$.

  At each step $n$ if $\tau$ is the $n$-the element of $\tr$
  we will additionally construct anchors 
  \begin{itemize}
  \item for each pair $X^n_\sigma$ and $X^n_{\pi(\sigma)}$
    such that $\sigma,\pi(\sigma)\leq\tau$
  \item and for each tripple
    $X^n_\sigma,X^n_{\pi(\sigma)},X^n_{\sigma^\smallfrown
      a}$ such that $a\in\omega+1$, $\pi(\sigma^\smallfrown
    a)\subseteq\pi(\sigma)$ and
    $\sigma,\pi(\sigma),\sigma^\smallfrown a\leq\tau$.
  \end{itemize}

  \subsection*{Completing the diagram.}

  As we now have a clear picture of what should be constructed
  let us argue that this is enough to finish the proof. For
  each $t\in\spc$ the intersection $\bigcap_n
  X_{t\upharpoonright n}$ has precisely one point so let us
  define $\varphi(t)$ to be this point. The other function,
  $\psi$ is defined as $f\circ\varphi\circ P^{-1}$. Let us
  check that this works. Both functions $\psi$ and $\varphi$
  are injective thanks to the disjoitness of the sets
  $C_\tau$ and to the fact that $X_\tau\subseteq
  f^{-1}[C_\tau]$. The function $\psi$ is open because
  $C_\tau$ are clopens.

  To see continuity of $\varphi$ notice first that since the
  sets $X_\tau$ have diameters vanishing to $0$, it suffices
  to check that $\varphi$ is continuous on each $\oplus^n$
  (which are treated as subsets of $\spc$ via the embedding
  $e: \tau\mapsto\tau^\smallfrown(\omega,\omega,\ldots)$).
  Continuity on $\oplus^n$ is checked inductively on the
  sets $S^n_k$ for $0\leq k\leq n$.

  The set $S^n_0$ consists of one point, so there is nothing
  to check. Suppose that $\tau_i\rightarrow \tau$,
  $\tau,\tau_i\in S^n_k,i\in\omega$. Then either the
  sequence is eventually constant or $\tau\in S^n_{k-1}$.
  Let us assume the latter.  By the inductive assumption and
  continuity of projection
  $\varphi(\pi(\tau_i))\rightarrow\varphi(\tau)$.  Now pick
  any $\varepsilon>0$. Let $m$ be such that
  $\diam(X_\sigma)<\varepsilon$ for $\sigma\in\oplus^m$ and
  $j\in\omega$ such that
  $d(\tau_j,\pi(\tau_j))<3^{-m}\varepsilon$. Let us write
  $\rho^\smallfrown\omega^l$ for $\rho$ extended by $l$ many
  $\omega$'s. By (\ref{ineq}) and coherence of projections
  we have
  $$h(X_{{\tau_j}^\smallfrown\omega^{m-n}},X_{\pi(\tau_j)^\smallfrown\omega^{m-n}})<\varepsilon,$$
  which implies that $\varphi(\tau_j)$ and
  $\varphi(\pi(\tau_j))$ are closer than $3\varepsilon$.
  This shows that $\varphi(\tau_j)\rightarrow\varphi(\tau)$
  and proves continuity of $\varphi$.

  \subsection*{Key lemma.}

  Now we state the key lemma, which will be used to
  guarantee ``continuity'' of the family of sets $X_\tau$.

  \begin{lemma}\label{limit}
    Let $X$ be a Borel set, $f:X\rightarrow \omega^\omega$ a
    Borel, not $\sigma$-continuous function. There exist
    a basic clopen $C_\omega\subseteq f[X]$
    % $\tau_\omega\in\omega^{<\omega}$ 
    and a compact set $X_\omega\subseteq
    f^{-1}[C_\omega]$ such that
    \begin{itemize}
    \item $f\har X_\omega$ is not $\sigma$-continuous,
    \item
      $X_\omega\subseteq\cl\big((f^{-1}[\omega^\omega\setminus
      C_\omega])^*\big)$.
    \end{itemize}
    The compact set $X_\omega$ can be chosen of arbitrarily
    small diameter.
  \end{lemma}
  \begin{proof} Without loss of generality assume that
    $f^{-1}[C]=(f^{-1}[C])^*$ for all clopen sets
    $C\subseteq\omega^\omega$.  Let us consider the
    following tree of open sets, indexed by $\otr$
    $$U_\tau=\int\big(f^{-1}[[\tau]]\big).$$ Let
    $G=\bigcap_n\bigcup_{|\tau|=n} U_\tau$ and
    $Z_\tau=f^{-1}[[\tau]]\setminus U_\tau$. Notice that
    $f\har G$ is continuous and since $X=G\cup\bigcup_\tau
    Z_\tau$ there is $\tau\in\otr$ such that $Z_\tau\not\in
    I_f$. Observe that
    $Z_\tau\subseteq\cl\big(\bigcup_{\tau'\not=\tau,|\tau'|=|\tau|}
    f^{-1}[[\tau']]\big)$ because if an open set $U\subseteq
    f^{-1}[[\tau]]$ is disjoint from
    $\bigcup_{\tau'\not=\tau,|\tau'|=|\tau|}
    f^{-1}[[\tau']]$ then $U\subseteq U_\tau$. Now put
    $C_\omega=[\tau]$ and pick any compact set with small
    diameter $X_\omega\subseteq Z_\tau$ such that
    $X_\omega\not\in I_f$.
  \end{proof}

  \subsection*{The construction.}

  We begin with $X_\emptyset=X^0_\emptyset=X$ and
  $C_\emptyset=\omega^\omega$. Without loss of generality
  assume that $X=X^*$. Suppose we have done $n-1$ steps of
  the inductive construction up $\tau\in\tr$. Let $|\tau|=l$
  and $\sigma=\tau\har(l-1)$. There are three cases.

  \textbf{Case 1.} The four points $\tau$, $\pi(\tau)$,
  $r(\tau)$ and $r(\pi(\tau))$ are equal. So
  $\tau=(\omega,\ldots,\omega)$ and $C_{\tau\upharpoonright
    n-1}$ and $X^{n-1}_\sigma$ are already constructed. In
  this case we use Lemma \ref{limit} to find a clopen set
  $C_\tau$ and a compact set $X_\tau\subseteq
  X^{n-1}_\sigma$ of diameter $<|\tau|\slash 3^{n+1}$ small
  enough so that no element of the anchors constructed so
  far is contained in $X_\tau$. We put $X^n_\tau=X_\tau^*$,
  $X^n_\sigma=(X^{n-1}_\sigma\setminus f^{-1}[C_\tau])^*$
  and $X^n_\rho=X^{n-1}_\rho$ for other $\rho<\tau$. By the
  assertion of Lemma \ref{limit} we still have
  $X^n_\tau\subseteq\cl(X^n_\sigma)$. In this case we do not
  need to construct any new anchors.

  \textbf{Case 2.} The two points $\pi(\tau)$ and $r(\tau)$
  are equal but distinct from $\tau$. Let
  $\delta=d(\tau,r(\tau))$. Since
  $X^{n-1}_{r(\tau)}\subseteq \cl (X^{n-1}_\sigma)$ by the
  inductive assumption, we may find finitely many sets
  $B_i\subseteq X^{n-1}_\sigma,i\leq k$ such that
  \begin{itemize}
  \item $h(\bigcup_i B'_i, X_{r(\tau)})<\delta$ for any
    $B'_i\subseteq B_i$,
  \item $B_i\not\in I_f$.
  \end{itemize}
  The second condition follows from
  $X^{n-1}_\sigma=(X^{n-1}_\sigma)^*$. We may assume that
  for each clopen set $C\subseteq 2^\omega$ the set $B_i\cap
  f^{-1}[C]$ is either empty or outside of the ideal $I_f$.

  We are going to find clopens $C_i\subseteq C_\sigma$, for
  $i\leq k$ such that $C_i\cap C_{r(\tau)}=\emptyset$ and
  then put $C_\tau=\bigcup_{i\leq k} C_i$,
  $X^n_\sigma=(X^{n-1}_\sigma\setminus\bigcup_i
  f^{-1}[C_i])^*$ and find $X_\tau\subseteq\bigcup_{i\leq k}
  B_i\cap f^{-1}[C_i]$.  We will have to carefully define
  $X^n_{r(\tau)}$ so that $X^n_{r(\tau)}\subseteq
  \cl(X^n_\sigma)$.

  It is easy to see that for any $A\subseteq X^{n-1}_\sigma$
  $$X^{n-1}_{r(\tau)}=X^{n-1}_{r(\tau)}\cap\cl\big((X^{n-1}_\sigma\cap A)^*\big)\ \cup\
  X^{n-1}_{r(\tau)}\cap\cl\big((X^{n-1}_\sigma\cap
  A^{c})^*\big)$$ so (putting
  $A=f^{-1}[C_\sigma\cap[(m,0)]]$ for $m<\omega$) we may
  inductively on $m$ pick binary sequences $\beta^m_i\in
  2^m,i\leq k$ such that $f^{-1}[[\beta^m_i]]\cap
  B_i\not=\emptyset$ and
  \begin{displaymath}
    X^{n-1}_{r(\tau)}\cap\cl\big((X^{n-1}_\sigma\setminus f^{-1}[\bigcup_{i\leq k} [\beta^m_i]])^*\big)\not\in I_f.
  \end{displaymath}

  We are going to carry on this construction up to some
  $m<\omega$ and put
  $X^n_\rho=\big(X^{n-1}_\rho\setminus\bigcup_{i\leq k}
  f^{-1}[[\beta^m_i]]\big)^*$ for
  $\rho<\tau,\rho\not\supseteq r(\tau)$ and
  $X^n_{\rho}=\big(X^{n-1}_{\rho}\cap\cl(X^n_\sigma)\big)^*$
  for $\rho<\tau,\rho\supseteq r(\tau)$.  We must, however,
  take care that this does not destroy the existing anchors.

  Since $f^{-1}[\{x\}]\in I_f$ for any $x\in 2^\omega$ and
  there are only finitely many elements of the existing
  anchors, we may pick $m<\omega$ and construct the
  sequences $\beta^m_i$ so that for any element $A$ of an
  anchor ``below'' $X^{n-1}_{r(\tau)}$ it is the case that
  $A\cap\cl(f^{-1}[C_\sigma\setminus\bigcup_{i\leq k}
  [\beta^m_i]]\cap X^{n-1}_\sigma)\not\in I_f$ and for any
  element $A$ of other anchors
  $A\setminus\big(\bigcup_{i\leq k}
  f^{-1}[\beta^m_i]\big)\not\in I_f$.

  Once we have constructed the sequences $\beta^m_i$ for
  $i\leq k$ we put $C_i=[\beta^m_i]$ and $C_\tau=\bigcup_i
  [\beta^m_i]$. Next we find $I_f$-positive compact sets
  $X_i$ inside $B_i\cap f^{-1}[[\beta_i]]$, each of diameter
  $<1\slash(3^{n+1}|\tau|)$.

  If $\delta>\diam(X^{n-1}_{\pi(\tau)})$ then we can pick
  one $X_i$ as $X_\tau$ and then
  $h(X_\tau,X^{n-1}_{\pi(\tau)})\leq 3\,
  h(X^{n-1}_\sigma,X^{n-1}_{\pi(\sigma)})<3^{|\tau|}\,\delta$.
  Otherwise, let $X_\tau=\bigcup_{i\leq k}X_i$ and then
  $\diam(X_\tau)<3\,\diam(X^{n-1}_{\pi(\tau)})\leq
  3\,\diam(X_{\pi(\tau)})$. Define $X^n_\tau=X_\tau^*$.

  At this step we create anchors for the pair $X^n_\tau$ and
  $X^n_{r(\tau)}$ as well as for the tripples $X^n_\sigma$,
  $X^n_\rho$, $X^n_\tau$ for $\rho<\tau$.

  \textbf{Case 3.} The two points $\pi(\tau)$, $r(\tau)$ are
  distinct. Let $\delta=d(\tau,\pi(\tau))$. By coherence of
  the projections $\pi(\tau)\supseteq\pi(\sigma)$. By the
  inductive assumption we have
  $h(X^{n-1}_\sigma,X^{n-1}_{\pi(\sigma)})<3^{|\sigma|}\,\delta$.
  Using the existing anchor for the tripple
  $X^{n-1}_\sigma$, $X^{n-1}_{\pi(\sigma)}$,
  $X^{n-1}_{\pi(\tau)}$ let us find finitely many sets
  $B_i,i\leq k$ in $X_\sigma$ such that
  \begin{itemize}
  \item $h(\bigcup_i B'_i,
    X_{\pi(\tau)})<3^{|\sigma|}\,\delta$ for any
    $B'_i\subseteq B_i$,
  \item $B_i\not\in I_f$.
  \end{itemize}
  As before, we assume assume that for each clopen set
  $C\subseteq 2^\omega$ if $B_i\cap f^{-1}[C]\in I_f$ then
  it is empty. We have now two subcases, in analogy to the
  two previous cases.

  \textbf{Subcase 3.1.} Suppose $\tau=r(\tau)$. Similarly as
  in Case 1, we use Lemma \ref{limit} to find $X_i\subseteq
  B_i$ and $C_i$ for $i\leq k$. Put $C_\tau=\bigcup_{i\leq
    k} C_i$. If $\delta>\diam(X^{n-1}_{\pi(\tau)})$ then we
  can pick one $X_i$ as $X_\tau$ and then
  $h(X_\tau,X^{n-1}_{\pi(\tau)})\leq 3\,
  h(X^{n-1}_\sigma,X^{n-1}_{\pi(\sigma)})<3^{|\tau|}\,\delta$.
  Otherwise, let $X_\tau=\bigcup_{i\leq k}X_i$ and then
  $\diam(X_\tau)<3\,\diam(X^{n-1}_{\pi(\tau)})\leq
  3\,\diam(X_{\pi(\tau)})$. Again, similarly as in Case 1,
  we put $X^n_\tau=(X^n_\tau)^*$,
  $X^n_\sigma=(X^{n-1}_\sigma\setminus f^{-1}[C_\tau])^*$,
  $X^n_\rho=X^{n-1}_\rho$ for other $\rho<\tau$.

  \textbf{Subcase 3.2.} Suppose $\tau\not=r(\tau)$.
  Similarly as in Case 2, we find clopens $C_i$ in
  $\omega^\omega$ such that $X^{n-1}_{r(\tau)}\cap
  \cl\big((f^{-1}[C_\sigma\setminus\bigcup_{i\leq
    k}C_i])^*\big)\not\in I_f$ and no existing anchor is
  destroyed when we put
  $X^n_\rho=(X^{n-1}_\rho\setminus\bigcup_{i\leq
    k}f^{-1}[C_i])^*$ for $\rho<\tau,\rho\not\supseteq
  r(\tau)$ and $X^n_\rho=X^{n-1}_\rho\cap\cl(X^n_\sigma)$
  for $\rho<\tau,\rho\supseteq r(\tau)$.

  Next we find $I_f$-positive compact sets $X_i\subseteq
  B_i\cap f^{-1}[C_i]$ each of diameter
  $<1\slash(3^{|\tau|+1}|\tau|)$. As previously, if
  $\delta>\diam(X^{n-1}_{\pi(\tau)})$ then we can pick one
  $X_i$ as $X_\tau$ and then
  $h(X_\tau,X^{n-1}_{\pi(\tau)})\leq 3\,
  h(X^{n-1}_\sigma,X^{n-1}_{\pi(\sigma)})<3^{|\tau|}\,\delta$.
  Otherwise, let $X_\tau=\bigcup_{i\leq k}X_i$ and then
  $\diam(X_\tau)<3\,\diam(X^{n-1}_{\pi(\tau)})\leq
  3\,\diam(X_{\pi(\tau)})$. Again, we put
  $X^n_\tau=(X^n_\tau)^*$,

  In Case 3 we construct the same anchors as in Case 2.

  This ends the construction and the entire proof.
\end{proof}

\begin{theorem}
  If $f:X\rightarrow\omega^\omega$ is not
  $\sigma$-continuous then there exist topological
  embeddings $\varphi$ and $\psi$ such that the following
  diagram commutes:
   $$
     \begin{CD}
       \omega^\omega @ >\psi >> \omega^\omega\\
       @AAP A @AAf A\\
       (\omega+1)^\omega @> \varphi >> X
     \end{CD}
   $$   
\end{theorem}
\begin{proof}
  By Theorem \ref{dichotomy} we have $\psi$ and $\varphi$
  such that $\psi$ is $1$-$1$ open. But as a Borel function
  it continuous on a dense $G_\delta$ set
  $G\subseteq\omega^\omega$. On the other hand by the
  properties of the function $P$ $X\in I_P$ implies $P[X]$
  is meager. So $P^{-1}[G]\not\in I_P$ and the problem
  reduces to the restriction of the function $P$.  This,
  however, has been proved in \cite{MS.for} (Corollary 2).
  So we get the following diagram:
   $$
     \begin{CD}
       \omega^\omega @>\psi'>> G @>\psi>> \omega^\omega\\
       @AAP A @AAP\upharpoonright G A @AAf A\\
       \spc@>\varphi'>> P^{-1}[G] @>\varphi>> X
     \end{CD}
   $$ 
   which ends the proof.
\end{proof}

\end{document}